\documentclass[11pt]{article}

\usepackage{amsbsy, graphicx, amssymb, amsmath, amsthm, amsfonts, mathrsfs, colortbl}
\usepackage[round]{natbib}
\usepackage[colorlinks, citecolor=blue, urlcolor=blue, linkcolor=blue, breaklinks]{hyperref}

\newtheorem{theorem}{Theorem}

\newtheorem{lemma}{Lemma}

\theoremstyle{remark}
\newtheorem{remark}{\bf Remark}

 \allowdisplaybreaks

\newcommand{\eps}{\varepsilon}

\newcommand{\md}{\mathrm{d}}

\begin{document}


\begin{center}
{\bf\LARGE  On the Identifiability Conditions in Some Nonlinear Time
Series Models }
 \vspace{.5cm}
\end{center}
\begin{center}
Jungsik Noh and Sangyeol Lee\\
Department of Statistics, Seoul National University\\
Email: nohjssunny@gmail.com; sylee@stats.snu.ac.kr
\end{center}

\begin{center}
    Revised February 28, 2015
\end{center}

\begin{abstract}
In this study, we consider the identifiability problem for nonlinear
time series models. Special attention is paid to smooth transition
GARCH, nonlinear Poisson autoregressive, and multiple regime smooth
transition autoregressive models. Some sufficient conditions are
obtained to establish the identifiability of these models.
\end{abstract}

\\
\noindent{\bf Abbreviated title}: Identifiability in nonlinear time
series models

\section{Introduction}\label{Introduction}

Verifying the identifiability conditions for time series models is a
fundamental task in constructing the consistent estimators of model
parameters and ensuring the positive definiteness of their
asymptotic covariance matrices. Although time series models are
assumed to be identifiable in many situations, its verification is
often nontrivial and even troublesome, especially in handling
nonlinear generalized autoregressive conditional heteroscedasticity
(GARCH) models. This issue has a long history and there exist a vast
amount of relevant studies in the literature. For instance,
\cite{Rothenberg1971} introduced the global and local identification
concept and verified that local identifiability is equivalent to the
nonsingularity of the information matrix. \cite{Phillips1989}
derived asymptotic theories in partially identified models.
\cite{Hansen1996} and \cite{FrancqHorvathZakoian2010} proposed a
test for the hypothesis wherein nuisance parameters are
unidentifiable. \cite{Komunjer2012} provided the primitive
conditions for global identification in moment restriction models.
In most cases, the identifiability condition is inherent to given
statistical models; for example, the multiple linear regression
model is unidentifiable when exact multicollinearity exists. Thus,
in nature, the verification of identifiability is more complicated
in nonlinear time series models with volatilities, such as threshold
autoregressive and smooth transition GARCH models (see, for
instance, \cite{Chan1993} and \cite{MeitzSaikkonen2011}). Thus,
there is a need to develop a more refined approach than the existing
ones to cope with the problem more adequately.

In this study, we deal with the identifiability problem within a
framework similar to that of the $M$-estimation. To elucidate, let
us consider the nonlinear least squares (NLS) estimation  from a
strictly stationary ergodic process $\{(Y_t, {{Z}}_t)\}$, with
$E(Y_t|{{Z}}_t)=f({{Z}}_t, \beta^\circ)$ for some known function
$f$. Then, the limit of the random objective functions for parameter
estimation is uniquely minimized at $\beta^\circ$ when the following
identifiability condition holds:
\begin{align}\label{NLS ID condition}
    f({Z}_1, \beta)=f({Z}_1, \beta^\circ) \text{ a.s. implies }  \beta=\beta^\circ.
\end{align}
In most $M$-estimation procedures, the identifiability conditions
are given in the form of (\ref{NLS ID condition}), where $f$ can be
a conditional mean, variance, or quantile function (see
\citet[p.~463]{Hayashi2000}, \cite{Berkesetal2003}, and
\cite{LeeNoh2013}).
Moreover, as seen in \cite{Wu1981}, to ensure the positive
definiteness of asymptotic covariance matrices of the NLS estimator,
one needs to verify that $\lambda^T
\partial{ f({Z}_1, \beta^\circ) }/\partial\beta=0$ a.s. implies
$\lambda=0$. The method described in this study is also useful to
verify the positive definiteness of asymptotic covariance matrices
of parameter estimators.

As a representative study on the issue with nonlinear time series,
we can refer to \cite{ChanTong1986}, who studied the asymptotic
theory of NLS estimators for the smooth transition AR (STAR) models
and verified the positive definiteness of asymptotic variance
matrices. Later, many authors handled this problem using various
GARCH-type models because it is crucial when verifying the
asymptotic properties of quasi-maximum likelihood estimators
(QMLEs). For example, \cite{StraumannMikosch2006},
\cite{MedeirosVeiga2009}, \cite{KristensenRahbek2009},
\cite{MeitzSaikkonen2011}, and \cite{LeeLee2012} consider the
identifiability problem in exponential and asymmetric GARCH($p,q$)
models, flexible coefficient GARCH($1,1$) models nesting  a smooth
transition GARCH($1,1$) (STGARCH) model, nonlinear ARCH models,
nonlinear AR($p$) models with nonlinear GARCH($1,1$) errors,
STAR($p$)-STGARCH($1,1$) models, and Box--Cox transformed threshold
GARCH($p,q$) models. To ensure (\ref{NLS ID condition}), these
authors developed their own methods that reflect the nonlinear
structure of underlying models.

In this study, we develop a method that refines existing ones to
deduce the identifiability conditions for various nonlinear time
series models, tribute to STGARCH($p,q$), Poisson autoregressive,
and multiple regime STAR($p$) models.
The remainder of this paper is organized as follows. In
Section~\ref{Sec2-heuristics}, we describe our method using some
examples. In Section~\ref{Sec3}, we investigate the identifiability
conditions in the aforementioned models. The proofs are provided in
the Appendix.

\section{Examples and motivation}\label{Sec2-heuristics}

In this section, we explore some existing methods that verify the
identifiability of STAR models and asymmetric GARCH (AGARCH) models.
 In what
follows, $\{X_t\}$ and $\mathcal{F}_t$ denote the data-generating
process and the $\sigma$-field generated by $\{X_s : s\leq t\}$.

First, we consider the STAR model with two regimes as follows:
\begin{align*}
    X_t &= m(X_{t-1}, \ldots, X_{t-p}; \theta^\circ)
    +\eps_t, \\
    m(X_{t-1}, \ldots, X_{t-p}; \theta^\circ) &= {\beta_0^\circ}^T \bold{X}_{t-1}  +  {\beta_1^\circ}^T \bold{X}_{t-1} F\left( {X_{t-d} - c^\circ \over z^\circ}  \right),
\end{align*}
where  $\{\eps_t\}$ are iid random variables, ${\theta^\circ}^T =
({\beta_0^\circ}^T, {\beta_1^\circ}^T, c^\circ, r^\circ)$ and
$\bold{X}_{t-1}=(1, X_{t-1}, \ldots, X_{t-p})^T$, and $F(\cdot)$ is
a smooth distribution function. \cite{ChanTong1986} verified the
positive definiteness of $E[\dot{m}_t(\theta^\circ)
\dot{m}_t(\theta^\circ)^T]$, where
$\dot{m}_t(\theta^\circ)=\dot{m}(X_{t-1}, \ldots, X_{t-p};
\theta^\circ)$ denotes the gradient of $m(x;\theta)$ at
$\theta^\circ$, by showing that for a given $\lambda\ne0$, there
exists $S\subset \mathbb{R}^p$, such that $\{ \lambda^T
\dot{m}(x;\theta^\circ)\}^2 $ is positive for any $x\in S$ and
$P(\left\{ (X_{t-1}, \ldots, X_{t-p}) \in S \right\})>0$. On the
other hand, \cite{MeitzSaikkonen2011} also considered the above
model  and verified that $m(X_{t-1}, \ldots, X_{t-p};
\theta)=m(X_{t-1}, \ldots, X_{t-p}; \theta^\circ)$ a.s. implies
$\theta=\theta^\circ$. In both cases, the main step is commonly to
show that the function $x \mapsto g(x;\theta, \theta^\circ)$, which
equals $(\theta-\theta^\circ)^T \dot{m}(x;\theta^\circ)$ in
\cite{ChanTong1986} and $m(x;\theta)- m(x;\theta^\circ)$ in
\cite{MeitzSaikkonen2011}, satisfies $g(x;\theta, \theta^\circ)=0$
for all $x\in supp(X_{t-1}, \ldots, X_{t-p})$, where $supp(Y)$
denotes the distribution support of the random vector $Y$. With this
equation, they could deduce certain conditions to guarantee
$\theta=\theta^\circ$. Motivated by these studies, we take a similar
approach to deduce the identifiability conditions for nonlinear time
series models. In fact, our method is handier than those in the
existing studies, such as \cite{KristensenRahbek2009},
\cite{MeitzSaikkonen2011}, and \cite{LeeLee2012}. For example, our
method no longer requires the condition that either the observations
or their conditional volatilities should take all values of an open
interval with a positive probability.

Next, we consider the case of an AGARCH($1$,$1$) model with power
$2$:
\begin{eqnarray}\label{AGARCH eqn}
    X_t = \sigma_t \eta_t,\ \   \sigma_t^2= \omega^\circ +  \alpha^\circ \left(|X_{t-1}|- \gamma^\circ X_{t-1} \right)^2
    + \beta^\circ \sigma_{t-1}^2,
\end{eqnarray}
where $\{ \eta_t \}$ is a sequence of iid random variables with
$E\eta_t=0$ and $E\eta_t^2=1$. \cite{KristensenRahbek2009} and
\cite{StraumannMikosch2006} derived identifiability conditions for
asymmetric power ARCH and AGARCH models. We denote $\theta^\circ = (
\omega^\circ, \alpha^\circ, \beta^\circ, \gamma^\circ)^T$ and
$\Theta = (0,\infty)\times [0,\infty) \times [0,1) \times [-1,1]$,
where $\alpha^\circ>0$. Assuming that Model~(\ref{AGARCH eqn}) has a
strictly stationary solution $\{X_t\}$,
for  $\theta\in\Theta$, we define a strictly stationary process
$\{\sigma_t^2(\theta)\}$ as the solution of
\begin{align}\label{def of sigma in AGARCH}
    \sigma_t^2(\theta) &=  \omega +  \alpha \left(|X_{t-1}|- \gamma X_{t-1} \right)^2
    +  \beta \sigma_{t-1}^2(\theta)
, \quad \forall t\in\mathbb{Z},
\end{align}
where $\sigma_t^2(\theta^\circ)$ is equal to $\sigma_t^2$.

In this case, the identifiability condition is that
$\sigma_t^2=\sigma_t^2(\theta)$ a.s. for some $t\in\mathbb{Z}$ and
$\theta\in\Theta$ implies $\theta=\theta^\circ$, which is crucial to
verify the strong consistency of QMLE. Below, we demonstrate the
approach of \cite{StraumannMikosch2006}. Note that
$\sigma_t^2=\sigma_t^2(\theta)$ a.s. for all $t$ because
$\{\sigma_t^2-\sigma_t^2(\theta)\}$ is stationary. Then, one can
obtain
\begin{align}\label{AGARCH equality1}
\begin{aligned}
    \omega^\circ - \omega + \sigma_{t-1}^2 Y_{t-1}  &= 0 \quad\text{a.s.},
\end{aligned}
\end{align}
where $Y_{t-1}= \alpha^\circ \left(|\eta_{t-1}|- \gamma^\circ
\eta_{t-1} \right)^2 - \alpha \left(|\eta_{t-1}|- \gamma \eta_{t-1}
\right)^2 + \beta^\circ - \beta$. As shown in Lemma~5.3 of
\cite{StraumannMikosch2006},  $Y_{t-1}$ is
$\mathcal{F}_{t-2}$-measurable due to (\ref{AGARCH equality1}), but
at the same time, it is independent of $\mathcal{F}_{t-2}$.
 Then,
$\theta=\theta^\circ$ can be easily deduced from the degeneracy of
$Y_{t-1}$ and certain mild conditions on the distribution of
$\eta_{t-1}$. This approach, however, cannot be extended
straightforwardly to more complicated models. Thus, in our study, we
take a different approach.

Our idea is to interpret the left-hand side of equation (\ref{AGARCH
equality1}) as a function of $\eta_{t-1}$. Considering that
$\sigma_{t-1}$ is given, for example, as constant $\sigma$,
we introduce the continuous function:
\begin{align*}
    g(x,\sigma) = \omega^\circ - \omega + \sigma^2 \left\{ \alpha^\circ \left(|x|- \gamma^\circ
x \right)^2 - \alpha \left(|x|- \gamma x \right)^2 + \beta^\circ -
\beta
    \right\}.
\end{align*}
Since (\ref{AGARCH equality1}) implies $g(\eta_{t-1},
\sigma_{t-1})=0$ a.s., it follows that $g(x,\sigma)=0$ for all
$(x,\sigma) \in supp(\eta_{t-1}, \sigma_{t-1} )$. Further, owing to
the independence of $\eta_{t-1}$ and $\sigma_{t-1}$, we have
$g(x,\sigma)=0$ for all $(x,\sigma) \in supp(\eta_{t-1})\times
supp(\sigma_{t-1} )$. This, in turn, implies
\begin{align}\label{deduction1}
    P\left\{ g(x, \sigma_{t-1} )=0  ~ \text{for all }x\in supp(\eta_{t-1}) \right\}=1.
\end{align}
Assume that $supp(\eta_{t-1})=\mathbb{R}$; in fact, it is sufficient
to assume that $supp(\eta_{t-1})$ comprises three distinct  (one
positive and one negative) real numbers. Then, $g(x, \sigma_{t-1}
)=0$ a.s. for all $x\in \mathbb{R}$ and, particularly $g(0,
\sigma_{t-1})= \omega^\circ - \omega + \sigma_{t-1}^2 \left(
\beta^\circ - \beta \right)=0$ a.s., which leads to
$\beta=\beta^\circ$ and $\omega=\omega^\circ$ owing to the
nondegeneracy of $\sigma_{t-1}^2$. Henceforth, the equation $g(x,
\sigma_{t-1} )=0$ a.s. $\forall x\in \mathbb{R}$  is now reduced to
\begin{align}\label{identity in sec2}
    \alpha^\circ \left(|x|- \gamma^\circ x \right)^2 - \alpha \left(|x|- \gamma x \right)^2=0, \quad\forall x\in\mathbb{R},
\end{align}
and thus, $\theta=\theta^\circ$ is  derived. This AGARCH($1,1$)
example demonstrates that equation (\ref{deduction1}) plays a
crucial role in obtaining the conditions to guarantee the
identifiability of a time series model. Later, to obtain the desired
results for general nonlinear time series models, such as STGARCH,
nonlinear Poisson autoregressive, and multiple regime STAR models,
we will often apply the equations analogous to (\ref{deduction1})
and  results such as $P\left\{ \lim_{x\to\infty} g(x, \sigma_{t-1}
)=0 \right\}=1$ or $P\left\{ \lim_{x\to-\infty} x^{-2} g(x,
\sigma_{t-1} )=0 \right\}=1$, as seen  in  the proof of
Theorem~\ref{STGARCH thm}.


\section{Identifiability in nonlinear time series}\label{Sec3}
\subsection{Smooth transition GARCH models}\label{Sec3.1-STGARCH}

\cite{Gonzalez-Rivera1998} introduced the STGARCH($p,q,d$) model:
\begin{align}\label{STGARCH eqn}
\begin{aligned}
    X_t &= \sigma_t \eta_t, \\
    \sigma_t^2 &= \omega^\circ + \sum_{i=1}^q \alpha^\circ_{1i} X_{t-i}^2 + \left( \sum_{i=1}^q \alpha^\circ_{2i} X_{t-i}^2 \right)
    F(X_{t-d}, \gamma^\circ) + \sum_{j=1}^p \beta^\circ_j \sigma_{t-j}^2,
\end{aligned}
\end{align}
where $\{\eta_t\}$ is the same as that in Model~(\ref{AGARCH eqn}),
\begin{align*}
    F(X_{t-d}, \gamma^\circ)= {1\over 1+ e^{\gamma^\circ X_{t-d}} } -
    {1\over2},
\end{align*}
$d\in\{1,\ldots,q\}$ is pre-specified, and $\gamma^\circ>0$ is the
smoothness parameter that determines the speed of transition. It is
noteworthy that when $\gamma^\circ\to\infty$, the STGARCH($1,1,1$)
model becomes a GJR-GARCH($1,1$) model proposed by
\cite{GlostenJagannathanRunkle1993}, which is identical to
Model~(\ref{AGARCH eqn}).

 We denote the true parameter vector by
$\theta^\circ=( \gamma^\circ, \omega^\circ, \alpha_{11}^\circ,
\ldots, \alpha_{1q}^\circ, \alpha_{21}^\circ, \ldots,
\alpha_{2q}^\circ, \beta_1^\circ, \ldots, \beta_p^\circ)^T$.  Let
 $\Theta=[0, \infty)\times (0,\infty)\times
A\times B $ be the parameter space, where
\begin{align}
    A&=\left\{ (\alpha_{11}, \ldots, \alpha_{1q}, \alpha_{21}, \ldots, \alpha_{2q})\in\mathbb{R}^{2q} :
\alpha_{1i}\geq0, ~|\alpha_{2i}|\leq2\alpha_{1i}, \forall i \right\}, \nonumber \\
    B&=\left\{ (\beta_1,\ldots, \beta_p)\in [0,1)^p : \sum_{j=1}^p \beta_j <1 \right\}, \label{sum beta<1 condition}
\end{align}
and assume that $\theta^\circ\in\Theta$ for the conditional variance
to be positive.

Sufficient conditions to ensure the existence of a stationary
solution for Model~(\ref{STGARCH eqn}) are not specified in the
literature. For instance, \cite{StraumannMikosch2006} and
\cite{MeitzSaikkonen2008} derived such conditions only for general
GARCH-type models.  However, for example, it can be seen that the
STGARCH($1,1,1$) model is stationary  when $E\left[ \log\left\{
\beta_1^\circ + \left(\alpha_{11}^\circ + {1 \over 2}
|\alpha_{21}^\circ| \right) \eta_{t-1}^2 \right\} \right]<0$ (cf.
Example~4 and Table~1 of \cite{MeitzSaikkonen2008}).

Given the stationary solution $\{X_t\}$ and a parameter vector
$\theta\in\Theta$, we define
\begin{align*}
    c_t(\alpha)=\omega + \sum_{i=1}^q \alpha_{1i} X_{t-i}^2 + \left( \sum_{i=1}^q \alpha_{2i} X_{t-i}^2 \right) F(X_{t-d}, \gamma),
\end{align*}
where $\alpha=(\gamma, \omega, \alpha_{11}, \ldots, \alpha_{1q},
\alpha_{21}, \ldots, \alpha_{2q} )$. Note that the polynomial
$\beta(z)=1-\sum_{j=1}^p \beta_j z^j$ has all its zeros outside the
unit disc because of (\ref{sum beta<1 condition}). Define
$\sigma_t^2(\theta) = \beta(B)^{-1} c_t(\alpha)$, where $B$ is the
backshift operator. Then, we have the following.

\begin{theorem}\label{STGARCH thm}
Let $\{X_t  \}$ be a  stationary process satisfying (\ref{STGARCH
eqn}) and suppose that
\begin{itemize}
  \item[(a)] $\alpha_{2i}^\circ \ne 0$ for some $1\leq i\leq q$ and $\gamma^\circ>0$.
  \item[(b)] The support of the distribution of $\eta_1$ is $\mathbb{R}$.
\end{itemize}
    Then, if $\sigma_t^2=\sigma_t^2(\theta)$ a.s. for some $t\in\mathbb{Z}$ and $\theta\in\Theta$, we have $\theta=\theta^\circ$.
\end{theorem}

\begin{remark}\label{remark1}
    It is remarkable that the identifiability in the STGARCH models needs no restriction
    concerning orders $p$ and $q$.
    The above theorem shows that the STGARCH($p,q,d$) models can be consistently estimated
    by fitting any STGARCH($p^*, q^*, d$) models with $p^*\geq p$ and $q^*\geq q$.
However, this is not true for GARCH and AGARCH models,
    wherein conditions such as (c) in Theorem~\ref{STGARCH thm partial iden} below are necessary.
    See   \cite{FrancqZakoian2004} and \cite{StraumannMikosch2006}.
\end{remark}

\begin{remark}
     As pointed out by a referee, the common root condition for the STGARCH models is not required owing to the reasons described below.
    Consider a STGARCH($0,1,d$) model and let $\sigma_t^2$ be
    the conditional variance. Multiplying $(1-\beta B)$ to both sides of the volatility equation,
    we get
    $ (1-\beta B) \sigma_t^2 = (1-\beta)\omega + \alpha_{11} X_{t-1}^2 - \beta \alpha_{11} X_{t-2}^2 + \alpha_{21}X_{t-1}^2 F(X_{t-d}, \gamma)
    -\beta \alpha_{21}X_{t-2}^2 F(X_{t-d-1}, \gamma)$.
    This, however, is not expressible as a form of STGARCH($1,2,d$) models, unlike we see in GARCH and AGARCH models.
\end{remark}

\begin{remark}
    As in the case of the AGARCH model in Section~\ref{Sec2-heuristics}, the support needs not be  $\mathbb{R}$. For example,
    $supp(\eta_1)=\mathbb{Z}$ is sufficient.
\end{remark}

Condition~(a) in Theorem~\ref{STGARCH thm} suggests that there
exists a smooth transition mechanism, that is, conditional variances
asymmetrically respond to positive and negative news. When it fails,
the STGARCH model becomes a standard GARCH model. The following
theorem demonstrates that model parameters in (\ref{STGARCH eqn})
are only partially identified when no such transition mechanism
exists.

\begin{theorem}\label{STGARCH thm partial iden}
    Let $\{X_t  \}$ be a  stationary process satisfying (\ref{STGARCH eqn}) with $\gamma^\circ=0$ or
    $\alpha_{2i}^\circ =0$, $i=1,\ldots,q$.
    Suppose that condition~(b) in Theorem~\ref{STGARCH thm} and the following condition hold:
    \begin{itemize}
        \item[(c)] $\alpha_{1i}^\circ >0$ for some $1\leq i\leq q$,  $(\alpha_{1q}^\circ, \beta_p^\circ) \ne (0,0)$, and
        the polynomials $\alpha_{1}^\circ(z)=\sum_{i=1}^q \alpha_{1i}^\circ z^i$
                and $\beta^\circ(z)=1-\sum_{j=1}^p \beta_j^\circ z^j$ have no common zeros.
    \end{itemize}
    If $\sigma_t^2=\sigma_t^2(\theta)$ a.s. for some $t\in\mathbb{Z}$ and $\theta\in\Theta$,
    then
    $\omega=\omega^\circ$, $\alpha_{1i}=\alpha_{1i}^\circ$, $\beta_j=\beta_j^\circ$
    for $1\leq i\leq q$, $1\leq j\leq p$, and either $\gamma=0$ or $\alpha_{2i}=0$, $1\leq i\leq q$ holds.
\end{theorem}

\begin{remark}
    The hypothesis testing of whether the smoothness mechanism exists has been studied by \cite{Gonzalez-Rivera1998}.
    This is a type of testing problem wherein nuisance parameters are unidentifiable under the null hypothesis.
    In addition, inference in a similar situation has been studied by
    \cite{Hansen1996} and \cite{FrancqHorvathZakoian2010}.
\end{remark}


\subsection{Threshold Poisson autoregressive models}\label{Sec3.2}

Poisson autoregressive models (or integer-valued GARCH models) are
used to model time series of counts with over-dispersion and have
been widely applied in fields ranging from finance to epidemiology
to estimate, for example, the number of transactions per minute of
certain stocks and the daily epileptic seizure counts of patients.
See \cite{FokianosRahbekTjostheim2009}, \cite{KangLee2014}, and the
references therein.

Let $\{X_t : t\geq0 \}$ be a time series of counts and $\{ \lambda_t
: t\geq0 \}$ its intensity process. Let $\mathcal{F}_{0,t}$ denote
the $\sigma$-field generated $\{ \lambda_0, X_0, \ldots, X_t \}$. An
integer-valued threshold GARCH (INTGARCH) model is then defined by
\begin{align}\label{INTGARCH eqn}
\begin{aligned}
    &X_t | \mathcal{F}_{0,t-1} \sim \text{Poisson}(\lambda_t),\\
    & \lambda_t =\omega^\circ + \alpha_1^\circ X_{t-1}+ (\alpha_2^\circ-\alpha_1^\circ)(X_{t-1} - l^\circ)^+ + \beta^\circ \lambda_{t-1}
\end{aligned}
\end{align}
for $t\geq1$, where $a^+$ denotes $\max\{0,a\}$. We assume that the
true parameter vector $\theta^\circ=(\omega^\circ, \alpha_1^\circ,
\alpha_2^\circ, \beta^\circ, l^\circ)$ belongs to a parameter space
$\Theta = (0,\infty)\times[0,1)^3\times\mathbb{N}$. Theorem~2.1 of
\cite{Neumann2011} indicates that if
$\beta^\circ+\max\{\alpha_1^\circ, \alpha_2^\circ\} <1$, there
exists a unique stationary bivariate process $\{ (X_t, \lambda_t) :
t\geq0 \}$ satisfying (\ref{INTGARCH eqn}). Then, the time domain
can be extended from $\mathbb{N}_0=\mathbb{N} \cup \{0\}$ to
$\mathbb{Z}$. \cite{FrankeKirchKamgaing2012} considered the
conditional LS estimation in these models.

Given the stationary process $\{X_t : t\in\mathbb{Z} \}$ and a
parameter vector $\theta\in\Theta$, we define a stationary process
$\{ \lambda_t(\theta) \}$ as the solution of
\begin{align*}
    \lambda_t(\theta) &=  \omega + \alpha_1  X_{t-1}
    + (\alpha_2 - \alpha_1 )  (X_{t-1}- l)^+ +\beta \lambda_{t-1}(\theta), \quad t\in\mathbb{Z}.
\end{align*}
Then, we have the following.

\begin{theorem}\label{INTGARCH thm}
  Suppose that $\{X_t: t\in\mathbb{Z} \}$ is a  stationary process satisfying (\ref{INTGARCH eqn}) and $\alpha_1^\circ\ne\alpha_2^\circ$.
  Then, if $\lambda_t=\lambda_t(\theta)$ a.s. for some $t\in\mathbb{Z}$ and $\theta\in\Theta$, we have $\theta=\theta^\circ$.
\end{theorem}

\begin{remark}
  When $\alpha_1^\circ=\alpha_2^\circ>0$, Model~(\ref{INTGARCH eqn}) becomes an integer-valued GARCH($1,1$) model.
  In this case, it can be seen that parameters, except the threshold parameter $l$, are identifiable.
\end{remark}

\subsection{General Poisson autoregressive models}\label{Sec3.3}

\cite{Neumann2011} considered a class of nonlinear Poisson
autoregressive models $\{X_t:t\in\mathbb{Z} \}$ of counts with
intensity process $\{\lambda_t:t\in\mathbb{Z} \}$ such as
\begin{align}\label{Poisson AR eqn}
\begin{aligned}
    &X_t | \mathcal{F}_{t-1} \sim \text{Poisson}(\lambda_t),  & & \lambda_t = f(\lambda_{t-1}, X_{t-1}, \theta^\circ )
\end{aligned}
\end{align}
for some known function $f:[0,\infty)\times \mathbb{N}_0 \times
\Theta \to [0,\infty)$. According to Theorems~2.1 and 3.1 of
\cite{Neumann2011}, when $f(\cdot,\theta^\circ)$ satisfies the
following contractive condition:
\begin{align*}
    |f(\lambda, y, \theta^\circ)-f(\lambda', y', \theta^\circ) |\leq \kappa_1|\lambda-\lambda'|+\kappa_2 |y-y'|,
    \quad\forall \lambda,\lambda'\geq0, ~\forall y,y'\in\mathbb{N}_0,
\end{align*}
where $\kappa_1, \kappa_2\geq0$ and $\kappa_1+\kappa_2 < 1$, there
exists a stationary process $\{(X_t,\lambda_t)\}$ with
$\lambda_t\in\mathcal{F}_{t-1}$  satisfying (\ref{Poisson AR eqn}).
Further, in view of Theorem~3.1 in \cite{Neumann2011}, one can
define a stationary process $\{\lambda_t(\theta)\}$ satisfying
\begin{align*}
    \lambda_t(\theta) = f(\lambda_{t-1}(\theta), X_{t-1}, \theta), \quad\forall t\in\mathbb{Z}
\end{align*}
for the stationary process $\{X_t\}$ and parameter vector
$\theta\in\Theta$. \cite{FokianosTjostheim2012} studied ML
estimation in these models.

The following theorem presents the mild requirements of $f$ for
their identifiability assumptions. Its proof is straightforward in
view of the proof of Theorem~\ref{INTGARCH thm}.

\begin{theorem}\label{Poisson AR thm}
  Let $\{(X_t, \lambda_t)\}$ be a stationary process satisfying (\ref{Poisson AR
  eqn}) and
  suppose that
\begin{itemize}
  \item[(a)] For each $\theta\in\Theta$, $f(\cdot,\theta)$ is continuous on $supp(\lambda_1)\times \mathbb{N}_0$.
  \item[(b)] $f(\lambda,y,\theta)=f(\lambda,y,\theta^\circ),~\forall\lambda\in supp(\lambda_1),~\forall y\in\mathbb{N}_0$
  implies $\theta=\theta^\circ$.
\end{itemize}
Then, if $\lambda_t=\lambda_t(\theta)$ a.s. for some
$t\in\mathbb{Z}$ and $\theta\in\Theta$, $\theta=\theta^\circ$.
\end{theorem}

\subsection{Multiple regime smooth transition autoregressive models}\label{Sec3.4}

Regime switching models for financial data have received
considerable attention. For example,  \cite{Terasvirta1994} studied
inference for two-regime STAR models and \cite{McAleerMedeiros2008}
and \cite{LiLing2012} considered multiple-regime smooth transition
and threshold AR models. In this subsection, we consider the
nonlinear LS estimation in a multiple-regime STAR model with
heteroscedastic errors proposed by \cite{McAleerMedeiros2008}.

Suppose that $\{X_t\}$ follows a multiple-regime STAR model of order
$p$ with $M+1$ (limiting) regimes, that is,
\begin{align}\label{STAR eqn}
    X_t = {\beta_0^\circ}^T \bold{X}_{t-1}  + \sum_{i=1}^M {\beta_i^\circ}^T \bold{X}_{t-1} G(X_{t-d^\circ} ; \gamma_i^\circ, c_i^\circ)
    +\eps_t,
\end{align}
where $\{\eps_t\}$ is white noise, $\beta_i^\circ=( \phi_{i0}^\circ,
\phi_{i1}^\circ, \ldots, \phi_{ip}^\circ )^T$
 for $0\leq i \leq M$,
$\bold{X}_{t-1}=(1, X_{t-1}, \ldots, X_{t-p})^T$, and
$G(X_{t-d^\circ} ; \gamma_i^\circ, c_i^\circ)$ is a logistic
transition function given by
\begin{align}\label{def of G()}
    G(X_{t-d^\circ} ; \gamma_i^\circ, c_i^\circ) &= \frac{1}{ 1+e^{-\gamma_i^\circ(X_{t-d^\circ} - c_i^\circ)}
    },
\end{align}
wherein the regime switches  according to the value of transition
variable $X_{t-d^\circ}$: $d^\circ \in \{1,\ldots,p\}$ is a delay
parameter, $-\infty<c_1^\circ<\cdots<c_M^\circ<\infty$ are threshold
parameters, and $\gamma_i^\circ>0$, $i=1,\ldots,M$, are  smoothing
parameters. When $\gamma_i^\circ$ is quite large, Model~(\ref{STAR
eqn}) is barely distinguishable from the threshold model studied by
\cite{LiLing2012}.

 In the literature, one can find sufficient
conditions under which Model~(\ref{STAR eqn})  is stationary when
the error terms are iid. For example, Theorem~2 of
\cite{McAleerMedeiros2008} ensures the stationarity of
Model~(\ref{STAR eqn}) of order 1. Using the same reasoning and
Lemma~2.1 of \cite{Berkesetal2003},
we can see that  
Model~(\ref{STAR eqn}) has a stationary solution if $$\sum_{j=1}^p
\sup_{x\in\mathbb{R} } \left| \phi^\circ_{0j} + \sum_{i=1}^M
\phi_{ij}^\circ G(x;\gamma^\circ_i, c^\circ_i) \right| < 1.$$ It is
also true if $\max_{0\leq i \leq M}  \sum_{j=1}^p \left|
\sum_{k=0}^i \phi^\circ_{kj}  \right| < 1$, which can be deduced
from Theorem~3.2 and Example~3.6 in \cite{AnHuang1996}.

We denote by $\theta=( \beta_0^T, \beta_1^T, \ldots, \beta_M^T,
\gamma_1, \ldots, \gamma_M, c_1, \ldots, c_M, d)^T$ a parameter
vector belonging to a parameter space $\Theta \subset
\mathbb{R}^{(M+1)(p+1)+2M}\times\{1,\ldots,p\}$ and set
\begin{align*}
    m(X_{t-1}, \ldots, X_{t-p}, \theta) &= {\beta_0}^T \bold{X}_{t-1}  + \sum_{i=1}^M {\beta_i}^T \bold{X}_{t-1} G(X_{t-d} ; \gamma_i, c_i).
\end{align*}
Then, we have the following.

\begin{theorem}\label{STAR thm}
  Let $\{X_t\}$ be a stationary process satisfying (\ref{STAR eqn}). Assume that
  \begin{itemize}
    \item[(a)] For each $i=1, \ldots, M$, $\beta_i^\circ \ne (0,\ldots, 0)^T \in\mathbb{R}^{p+1}$.
    \item[(b)] The support of the stationary distribution of $(X_p, \ldots, X_1)$ is $\mathbb{R}^p$.
    \item[(c)] The parameter space $\Theta$ satisfies that $\gamma_i>0$, $i=1,\ldots,M$, and $-\infty<c_1<\cdots<c_M<\infty$.
  \end{itemize}
Then, if $m(X_{t-1}, \ldots, X_{t-p}, \theta^\circ)=m(X_{t-1},
\ldots, X_{t-p}, \theta)$ a.s. for some $t\in\mathbb{Z}$ and
$\theta\in\Theta$, we have $\theta=\theta^\circ$.
\end{theorem}

\begin{remark}
  Theorem~\ref{STAR thm} is closely related to the identifiability of the finite mixture of logistic distributions
  (see Lemma~\ref{Linear independence Lemma} in the Appendix).
  Although the restriction on threshold parameters has a natural interpretation,
  it is not necessarily required. In fact, if we only assume that $(\gamma_i^\circ, c_i^\circ), i=1,\ldots,M$, are distinct,  
  instead of the condition $c_i^\circ<c_{i+1}^\circ$,
  then Model~(\ref{STAR eqn}) is weakly identifiable in the sense of \cite{RednerWalker1984}.
\end{remark}

\appendix

\section*{Appendix}\label{appendix}
\setcounter{equation}{0}

\renewcommand{\theequation}{A.\arabic{equation}}
\renewcommand{\thelemma}{A.\arabic{lemma}}
\renewcommand{\thecorollary}{A.\arabic{corollary}}

\textbf{Proof of Theorem~\ref{STGARCH thm}. } We only prove the
theorem when $d=1$ since
the other cases can be handled similarly. Owing to the stationarity,
we have $\sigma_t^2=\sigma_t^2(\theta)$ a.s. for any
$t\in\mathbb{Z}$. Since $\beta^\circ(z)\ne0$ for $|z|\leq1$ and
$\sigma_t^2=\beta^\circ(B)^{-1} c_t(\alpha^\circ)$, we can express
\begin{align}\label{STGARCH equality1}
    c_t(\alpha) &= \beta(B) \sigma_t^2(\theta)= \beta(B) \beta^\circ(B)^{-1} c_t(\alpha^\circ) = c_t(\alpha^\circ) + \sum_{j=1}^\infty b_j c_{t-j}(\alpha^\circ),
\end{align}
where $1+\sum_{j=1}^\infty b_j z^j = \beta(z)/\beta^\circ(z)$ for
$|z|\leq 1$. As discussed in Section~\ref{Sec2-heuristics}, we can
express (\ref{STGARCH equality1}) as a function of $\eta_{t-1}$ and
$\mathcal{F}_{t-2}$-measurable random variables:
\begin{align*}
\begin{aligned}
    &g_1( \eta_{t-1}, \sigma_{t-1}, A_{t,2}, B_{t,2}, A_{t,2}^\circ, B_{t,2}^\circ, D_{t,2} ) \\
    &:= (\alpha_{11}-\alpha_{11}^\circ) \sigma_{t-1}^2 \eta_{t-1}^2 + A_{t,2}-A_{t,2}^\circ
    + \left( \alpha_{21}\sigma_{t-1}^2 \eta_{t-1}^2 + B_{t,2} \right)
    F(\sigma_{t-1}\eta_{t-1}, \gamma) \\
    &\quad - \left( \alpha_{21}^\circ \sigma_{t-1}^2 \eta_{t-1}^2 + B_{t,2}^\circ \right)
    F(\sigma_{t-1}\eta_{t-1}, \gamma^\circ) - D_{t,2} \\
    &= 0 \quad{\text a.s.},
\end{aligned}
\end{align*}
where for $2\leq i^* \leq q$ and $2\leq k$,
\begin{align*}
    A_{t,i^*} &= \omega + \sum_{i=i^*}^q \alpha_{1i} X_{t-i}^2, &
    B_{t,i^*} &= \sum_{i=i^*}^q \alpha_{2i} X_{t-i}^2, &
    D_{t,k} &= \sum_{j=k-1}^\infty b_j c_{t-j}(\alpha^\circ), \\
    A_{t,i^*}^\circ &= \omega^\circ + \sum_{i=i^*}^q \alpha_{1i}^\circ X_{t-i}^2, &
    B_{t,i^*}^\circ &= \sum_{i=i^*}^q \alpha_{2i}^\circ X_{t-i}^2.
\end{align*}
Using the arguments that obtain (\ref{deduction1}) and
condition~(b), we can see that  with probability $1$,  $
    g_1(x, \sigma_{t-1}, A_{t,2}, B_{t,2}, A_{t,2}^\circ, B_{t,2}^\circ, D_{t,2} )=0
$ for all $x\in\mathbb{R}$. Particularly, this implies
\begin{align}\label{g_1(0)=0}
    g_1(0, \sigma_{t-1}, A_{t,2}, B_{t,2}, A_{t,2}^\circ, B_{t,2}^\circ, D_{t,2} )=A_{t,2}-A_{t,2}^\circ-D_{t,2}=0 \quad\text{a.s..}
\end{align}
Then, viewing (\ref{g_1(0)=0}) as a function of $\eta_{t-2}$ and
$\mathcal{F}_{t-3}$-measurable random variables, we can express
\begin{align}\label{STG-g_2(x)=0}
\begin{aligned}
    &g_2(\eta_{t-2}, \sigma_{t-2}, A_{t,3}, A_{t,3}^\circ, A_{t-1,2}^\circ, B_{t-1,2}^\circ, D_{t,3} ) \\
    &:= ( \alpha_{12} - \alpha_{12}^\circ ) \sigma_{t-2}^2 \eta_{t-2}^2 + A_{t,3}-A_{t,3}^\circ
    - b_1 c_{t-1}(\alpha^\circ) - D_{t,3} \\
    &= ( \alpha_{12} - \alpha_{12}^\circ ) \sigma_{t-2}^2 \eta_{t-2}^2 + A_{t,3}-A_{t,3}^\circ - D_{t,3} \\
    &\quad - b_1 \left\{ \alpha_{11}^\circ \sigma_{t-2}^2 \eta_{t-2}^2 + A_{t-1,2}^\circ + \left( \alpha_{21}^\circ \sigma_{t-2}^2 \eta_{t-2}^2 + B_{t-1,2}^\circ \right) F(\sigma_{t-2}\eta_{t-2}, \gamma^\circ)
    \right\}
     \\
    &=0 \quad\text{a.s.},
\end{aligned}
\end{align}
which entails
\begin{align}\label{STG-prob of iden is one}
    P\left( g_2(x, \sigma_{t-2}, A_{t,3}, A_{t,3}^\circ, A_{t-1,2}^\circ, B_{t-1,2}^\circ, D_{t,3} )=0, \forall x\in\mathbb{R}
    \right)=1.
\end{align}
Note that if
\begin{align}\label{STG-identity}
    f(x):= ax^2 + b + (cx^2 + d) F(\sigma x, \gamma^\circ)=0
\end{align}
for all $x\in\mathbb{R}$, where $a, b, c, d, \sigma>0,
\gamma^\circ>0$ are real numbers, because $\lim_{x\to\pm\infty}
x^{-2}f(x)=0$ and $\lim_{x\to\pm\infty}f(x)=0$, it must hold that
$a=c=0$ and $b=d=0$. Then, combining this and (\ref{STG-prob of iden
is one}), we get $b_1 \alpha_{21}^\circ=0$ and $b_1
B_{t-1,2}^\circ=0$ a.s.. Further, $B_{t-1,2}^\circ=0$ a.s. if and
only if $\alpha_{22}^\circ=\cdots=\alpha_{2q}^\circ=0$. Due to
condition~(a) and (\ref{STG-g_2(x)=0}), we have $b_1=0$ and
$A_{t,3}-A_{t,3}^\circ-D_{t,3}=0$ a.s., and similarly, it can be
seen that $b_k=0$, $k\geq2$, $A_{t,k+2}-A_{t,k+2}^\circ-D_{t,k+2}=0$
a.s., $2\leq k \leq q-2$, and $\omega-\omega^\circ-D_{t,k+2}=0$
a.s., $k\geq q-1$. This implies $\beta(\cdot) = \beta^\circ(\cdot)$,
$\omega=\omega^\circ$, and $A_{t,2}=A_{t,2}^\circ, \ldots,
A_{t,q}=A_{t,q}^\circ$ a.s., and subsequently,
$\alpha_{1q}=\alpha_{1q}^\circ, \ldots,
\alpha_{12}=\alpha_{12}^\circ$. From this and (\ref{STGARCH
equality1}), we can obtain
\begin{align}\label{STG-h_1(x)}
\begin{aligned}
    h_1( \eta_{t-1}, \sigma_{t-1},  B_{t,2},  B_{t,2}^\circ )
    &:= (\alpha_{11}-\alpha_{11}^\circ) \sigma_{t-1}^2 \eta_{t-1}^2
    + \left( \alpha_{21}\sigma_{t-1}^2 \eta_{t-1}^2 + B_{t,2} \right)
    F(\sigma_{t-1}\eta_{t-1}, \gamma) \\
    &\quad - \left( \alpha_{21}^\circ \sigma_{t-1}^2 \eta_{t-1}^2 + B_{t,2}^\circ \right)
    F(\sigma_{t-1}\eta_{t-1}, \gamma^\circ)  \\
    &= 0 \quad{\text a.s.}
\end{aligned}
\end{align}

Suppose that $\gamma=0$. Then, $F(X_{t-1},\gamma)\equiv0$, and using
 (\ref{STG-identity}) and (\ref{STG-h_1(x)}), we get
$\alpha_{21}^\circ=0$ and $B_{t,2}^\circ=0$ a.s. Since this is a
contradiction to  condition~(a),  $\gamma$ must  be positive. Thus,
from (\ref{STG-h_1(x)}), we have
\begin{align*}
    \lim_{x\to\infty} x^{-2} h_1(x, \sigma_{t-1},  B_{t,2},  B_{t,2}^\circ )
    = \sigma_{t-1}^2 \left\{ \alpha_{11} - \alpha_{11}^\circ - 2^{-1} (\alpha_{21}- \alpha_{21}^\circ)
    \right\}=0
    \quad\text{a.s..}
\end{align*}
Further, taking the limit $x\to -\infty$, we obtain
$\alpha_{11}=\alpha_{11}^\circ$ and $\alpha_{21}=\alpha_{21}^\circ$,
so that
\begin{align*}
    \lim_{x\to\infty}  h_1(x, \sigma_{t-1},  B_{t,2},  B_{t,2}^\circ )
    = -2^{-1} B_{t,2}+2^{-1} B_{t,2}^\circ=0 \quad\text{a.s.},
\end{align*}
which results in $\alpha_{2i}=\alpha_{2i}^\circ,~2\leq i\leq q$.
Then, in view of (\ref{STG-h_1(x)}), we obtain
\begin{align*}
    h_2(\eta_{t-1}, \sigma_{t-1}, B_{t,2}^\circ)
    &:= \left( \alpha_{21}^\circ \sigma_{t-1}^2 \eta_{t-1}^2 + B_{t,2}^\circ \right)
    \left( \frac{1}{1+e^{\gamma \sigma_{t-1}\eta_{t-1} } } - \frac{1}{1+e^{\gamma^\circ \sigma_{t-1}\eta_{t-1} } }
    \right)=0 \quad\text{a.s..}
\end{align*}
If $\gamma<\gamma^\circ$ and additionally if
$\alpha_{21}^\circ\ne0$, we should have
\begin{align*}
    \lim_{x\to\infty} x^{-2} e^{\gamma \sigma_{t-1} x  } h_2(x, \sigma_{t-1}, B_{t,2}^\circ)
    =  \alpha_{21}^\circ \sigma_{t-1}^2  =0
    \quad\text{a.s.},
\end{align*}
which leads to a contradiction. However, if $\alpha_{21}^\circ=0$,
we have $\lim_{x\to\infty} e^{\gamma \sigma_{t-1} x  } h_2(x,
\sigma_{t-1}, B_{t,2}^\circ)  = B_{t,2}^\circ=0$ a.s., which also
leads to a contradiction to condition~(a). Hence, we must have
$\gamma\geq\gamma^\circ$. Since $\gamma>\gamma^\circ$ is also
impossible, we conclude that $\gamma=\gamma^\circ$, which completes
the proof.
\hfill \fbox{} \\

\textbf{Proof of Theorem~\ref{STGARCH thm partial iden}. } As in
handling (\ref{STG-g_2(x)=0}), we follow the same lines in the proof
of Theorem~\ref{STGARCH thm} to obtain
\begin{align*}
\begin{aligned}
    &g'_2(\eta_{t-2}, \sigma_{t-2}, A_{t,3}, A_{t,3}^\circ, A_{t-1,2}^\circ,  D_{t,3} ) \\
    &:= ( \alpha_{12} - \alpha_{12}^\circ - b_1 \alpha_{11}^\circ  ) \sigma_{t-2}^2 \eta_{t-2}^2
    + A_{t,3}-A_{t,3}^\circ - b_1 A_{t-1,2}^\circ - D_{t,3} \\
    &=0 \quad\text{a.s..}
\end{aligned}
\end{align*}
Then, as in handling (\ref{g_1(0)=0}), we get $A_{t,3}-A_{t,3}^\circ
- b_1 A_{t-1,2}^\circ - D_{t,3}=0$ a.s. Similarly, it can be seen
that $ \omega-\omega^\circ -  b_1 A_{t-1,q}^\circ - b_2
A_{t-2,q-1}^\circ - \cdots - b_{q-1} A_{t-q+1,2}^\circ - D_{t,q+1}=0
$ a.s. Then, with probability $1$, for all $x\in\mathbb{R}$,
\begin{align}\label{STG g(x)=0}
\begin{aligned}
    g(x) &:= (\omega-\omega^\circ)
    - b_1 \left( \alpha_{1q}^\circ \sigma_{t-q-1}^2 x^2 + \omega^\circ \right)
    - b_2 \left( \alpha_{1,q-1}^\circ \sigma_{t-q-1}^2 x^2 + A_{t-2,q}^\circ \right)
    -\cdots \\
    &\quad - b_{q-1} \left( \alpha_{12}^\circ \sigma_{t-q-1}^2 x^2 + A_{t-q+1,3}^\circ  \right)
     - b_q \left( \alpha_{11}^\circ \sigma_{t-q-1}^2 x^2 + A_{t-q,2}^\circ    \right)
     - D_{t,q+2} \\
    &=0,
\end{aligned}
\end{align}
which, in turn, implies
\begin{align*}
    P\left( \lim_{x\to\infty}{-g(x) \over \sigma_{t-q-1}^2 x^2 }
    = b_1 \alpha_{1q}^\circ + \cdots + b_q  \alpha_{11}^\circ
    =0
    \right)=1.
\end{align*}
In fact,  we can obtain an analogous relationship between
$\eta_{t-q-k}$ and $\mathcal{F}_{t-q-k-1}$-measurable random
variables, $k\geq2$, and as such, $b_k \alpha_{1 q}^\circ + \cdots +
b_{k+q-1} \alpha_{11}^\circ=0$ for all $k\geq1$, which implies that
$
     \beta(z)  \beta^\circ(z)^{-1}  \alpha_1^\circ(z)
$ is a polynomial of  at most $q$ orders. Then, using condition~(c)
and the arguments similar to those in Straumann and Mikosch~(2006),
p. 2481, we  can see that $\beta(\cdot) = \beta^\circ(\cdot)$, and
thus, $b_j=0$ for $j\geq1$. Combining this, (\ref{g_1(0)=0}) and
(\ref{STG g(x)=0}), we get $A_{t,2}=A_{t,2}^\circ, \ldots,
A_{t,q}=A_{t,q}^\circ$ a.s. and $\omega=\omega^\circ$, which, in
turn, implies $\alpha_{1q}=\alpha_{1q}^\circ, \ldots,
\alpha_{12}=\alpha_{12}^\circ$. Hence, (\ref{STGARCH equality1}) can
be reexpressed as
\begin{align*}
\begin{aligned}
    h'_1( \eta_{t-1}, \sigma_{t-1},  B_{t,2} )
    &:= (\alpha_{11}-\alpha_{11}^\circ) \sigma_{t-1}^2 \eta_{t-1}^2
    + \left( \alpha_{21}\sigma_{t-1}^2 \eta_{t-1}^2 + B_{t,2} \right)
    F(\sigma_{t-1}\eta_{t-1}, \gamma) \\
    &= 0 \quad{\text a.s.}.
\end{aligned}
\end{align*}
From this, we can easily obtain $\alpha_{11}=\alpha_{11}^\circ$ and
the same equation as in (\ref{STG-identity}), which finally leads to
$\alpha_{21}=\cdots=\alpha_{2q}=0$. This completes the proof.
\hfill \fbox{} \\

\textbf{Proof of Theorem~\ref{INTGARCH thm}. } First, we conjecture
that the support of the stationary distribution of $(X_1,
\lambda_1)$ is a Cartesian product of $\mathbb{N}_0$ and
$supp(\lambda_1)$. If it is not true, there exists $(m',
\lambda')\in\mathbb{N}_0 \times supp(\lambda_1)$ such that
$(m',\lambda') \notin supp(X_1,\lambda_1)$, and for some positive
real number $r$,
\begin{align*}
    0 = P\left(  X_1=m', \lambda_1 \in (\lambda'-r,\lambda'+r)
    \right)
    = \int_{\lambda'-r}^{\lambda'+r} (m'!)^{-1} e^{-u} u^{m'} \md F_{\lambda_1} (u),
\end{align*}
where $F_{\lambda_1}$ is the distribution function of $\lambda_1$.
Since the integrand is positive, it must hold that $P(\lambda_1 \in
(\lambda'-r,\lambda'+r) )=0$, which, however, contradicts to the
fact that $\lambda'\in supp(\lambda_1)$. Thus, our conjecture is
validated.

Note that owing to the stationarity, for all $t\in\mathbb{Z}$,
\begin{align*}
\begin{aligned}
    &g(X_{t-1}, \lambda_{t-1})\\
    &:= (\omega-\omega^\circ) + (\alpha_1-\alpha_1^\circ) X_{t-1} + (\alpha_2 - \alpha_1 )  (X_{t-1}- l)^+
    - (\alpha_2^\circ - \alpha_1^\circ )  (X_{t-1}- l^\circ)^+
    + (\beta-\beta^\circ) \lambda_{t-1}\\
    &=0 \quad\text{a.s.,}
\end{aligned}
\end{align*}
and therefore,
\begin{align}\label{INTGARCH equality1}
    g(m,\lambda)=0 \quad\text{for all }m\in\mathbb{N}_0 \text{ and }\lambda\in
    supp(\lambda_1),
\end{align}
since $g(\cdot)$ is continuous and $supp(X_{1},
\lambda_{1})=\mathbb{N}_0\times supp(\lambda_1)$. In particular,
$g(0,\lambda)=(\omega-\omega^\circ)+(\beta-\beta^\circ) \lambda=0$
for any $\lambda\in supp(\lambda_1)$. Note that $\lambda_t$ is not
degenerate when $\alpha_1^\circ\ne \alpha_2^\circ$, since otherwise,
$X_{t-1}$ should be degenerate. Thus, we have $\omega=\omega^\circ$
and $\beta=\beta^\circ$, so that
$g(1,\lambda)=\alpha_1-\alpha_1^\circ=0$. Further, it follows from
(\ref{INTGARCH equality1}) that $\lim_{m\to\infty} m^{-1}
g(m,\lambda)=\alpha_2-\alpha_2^\circ=0$. Then, using the fact that
$g(l,\lambda)=g(l^\circ,\lambda)=0$ and $\alpha_1^\circ\ne
\alpha_2^\circ$, we obtain $l=l^\circ$, which completes the proof.
\hfill \fbox{} \\

\textbf{Proof of Theorem~\ref{STAR thm}. } For simplicity, we assume
that $d^\circ=1$:  the other cases can be handled similarly. From
condition (b) and the continuity of $m(\cdot, \theta)$, we can see
that
\begin{align}\label{STAR equality1}
    m(x_1, \ldots, x_p, \theta^\circ) = m(x_1, \ldots, x_p, \theta), \quad\forall x_j \in\mathbb{R},~1\leq j \leq p.
\end{align}
Suppose that $d\ne 1$.  From (\ref{STAR equality1}), we can express
\begin{align}\label{STAR d not d0 equality}
\begin{aligned}
    &m(x_1, \ldots, x_p, \theta^\circ) - m(x_1, \ldots, x_p, \theta) \\
    &= \left\{ f_0^\circ(\mathbf{x}_2) - f_0(\mathbf{x}_2) - \sum_{i=1}^M f_i(\mathbf{x}_2) G(x_d;\gamma_i, c_i)  \right\} \\
    &\quad+ \left\{ \phi_{0 1}^\circ - \phi_{0 1} - \sum_{i=1}^M \phi_{i 1} G(x_d; \gamma_i, c_i) \right\} x_{1}
     + \sum_{i=1}^M  \left( f_i^\circ(\mathbf{x}_2) + \phi_{i 1}^\circ x_{1} \right) G(x_{1};\gamma_i^\circ, c_i^\circ) \\
    &= 0,
\end{aligned}
\end{align}
where $G(\cdot)$ is the one in (\ref{def of G()}), $\mathbf{x}_2=(x_2, \ldots, x_p)^T$, and  
\begin{align*}
    f_i^\circ(\mathbf{x}_2) &= \phi_{i0}^\circ+ \sum_{2\leq j \leq p} \phi_{ij}^\circ x_j, &
    f_i(\mathbf{x}_2) &= \phi_{i0}+ \sum_{2\leq j \leq p} \phi_{ij} x_j, & &\text{for }i=0, 1, \ldots, M.
\end{align*}
Then, applying Lemma~\ref{Linear independence Lemma} below to
(\ref{STAR d not d0 equality}), we have $\phi_{1 1}^\circ=0$ and
$f_1^\circ(\mathbf{x}_2)=0$ for each
$\mathbf{x}_2\in\mathbb{R}^{p-1}$, which, however, contradicts to
condition~(a). Thus, it must hold that $d=d^\circ=1$. Owing to the
above, we can reexpress (\ref{STAR equality1}) as
\begin{multline}\label{STAR equality2}
    \left( f_0^\circ(\mathbf{x}_2) +\phi_{0 1}^\circ x_{1} \right)+ \sum_{i=1}^M  \left( f_i^\circ(\mathbf{x}_2) + \phi_{i 1}^\circ x_{1} \right) G(x_{1};\gamma_i^\circ, c_i^\circ) \\
    =
    \left( f_0(\mathbf{x}_2) + \phi_{0 1} x_{1} \right)+ \sum_{i=1}^M  \left( f_i(\mathbf{x}_2) + \phi_{i 1} x_{1} \right) G(x_{1};\gamma_i, c_i)
    , \quad\forall x_j \in\mathbb{R},~1\leq j \leq p.
\end{multline}
Lemma~\ref{Linear independence Lemma} ensures that  
a family of real-valued functions $    \mathcal{G}=\{ 1, i(\cdot) \}
\cup \{ G(\cdot;\gamma, c) : \gamma>0, c\in\mathbb{R} \} \cup \{
i(\cdot) G(\cdot;\gamma, c) : \gamma>0, c\in\mathbb{R}
    \},
$ where $i(\cdot)$ is an identity function, i.e., $i(y)=y$, are
linearly independent. Thus, any element of the linear span of
$\mathcal{G}$ is uniquely represented as a linear combination of the
elements of $\mathcal{G}$: see \cite{YakowitzSpragins1968}. Further,
there exists a vector $\mathbf{x}'_2 \in\mathbb{R}^{p-1}$ such that
$(f_i^\circ( \mathbf{x}'_2 ), \phi_{i1}^\circ) \ne (0,0)$ for all
$i=1,\ldots, M$; unless otherwise,
$\phi_{i0}^\circ=\cdots=\phi_{ip}^\circ=0$ for some $i$, which
contradicts condition~(a). Then, viewing (\ref{STAR equality2}) with
$\mathbf{x}_2$ substituted by $\mathbf{x}'_2$ as a function of $x_1$
and using condition~(c), we obtain $\phi_{01}^\circ=\phi_{01}$ and
$\phi_{i1}^\circ=\phi_{i1}$, $\gamma_i^\circ=\gamma_i$,
$c_i^\circ=c_i$ for $i=1,\ldots,M$. Subsequently, owing to
(\ref{STAR equality2}),  for all $x_1\in\mathbb{R}$ and
$\mathbf{x}_2\in\mathbb{R}^{p-1}$, we get
\begin{align*}
    \left( f_0^\circ(\mathbf{x}_2) - f_0(\mathbf{x}_2) \right) + \sum_{i=1}^M  \left( f_i^\circ(\mathbf{x}_2) - f_i(\mathbf{x}_2) \right)
    G(x_{1};\gamma_i^\circ, c_i^\circ)
    = 0.
\end{align*}
Then, applying Lemma~\ref{Linear independence Lemma} again, we
conclude that $\phi_{i0}^\circ=\phi_{i0}$ and
$\phi_{ij}^\circ=\phi_{ij}$,  $j=2,\ldots,p$, $i=0, 1,\ldots,M$.
This completes the proof.
\hfill \fbox{} \\

\begin{lemma}\label{Linear independence Lemma}
  Let $(\gamma_1, c_1), \ldots, (\gamma_k, c_k)$ be distinct real vectors with $\gamma_i>0$, $i=1,\ldots,k$.
  Suppose that for all $y\in\mathbb{R}$,
  \begin{align}\label{Lemma equality}
    d_{00}+d_{01}y + \sum_{i=1}^k (d_{i0} + d_{i1}y ) \frac{1}{ 1+e^{-\gamma_i (y-c_i)} } = 0.
  \end{align}
  Then, $d_{i0}=d_{i1}=0$ for $i=0, 1, \ldots, k$.
\end{lemma}

\textbf{Proof.} Denote by $g(y)$ the left-hand side of (\ref{Lemma
equality}). Then,  $\lim_{y\to -\infty} y^{-1}g(y)=d_{01}=0$, and
thus, $\lim_{y\to -\infty} g(y)=d_{00}=0$.  In what follows, for
function $f:\mathbb{R}\to\mathbb{R}$, we denote by
$\mathcal{L}\{{f}\}$ its two-sided Laplace transform, that is,
$\mathcal{L}\{f(\cdot)\}(s)=\int_{-\infty}^{\infty} e^{-sy} f(y) \md
y$. Note that the transform of the logistic distribution function is
as follows:
\begin{align*}
    F_0(s;\gamma,c) :=
    \mathcal{L}\{G(\cdot;\gamma, c)\}(s) = \frac{\pi \gamma^{-1} e^{-cs} }{\sin{ \pi \gamma^{-1}s} }, \quad 0<s< \gamma.
\end{align*}
Further,
\begin{align*}
    F_1(s;\gamma,c) :=
    \mathcal{L}\{i(\cdot) G(\cdot ;\gamma, c)\}(s)
    = {\pi \gamma^{-1} c e^{-cs} \over \sin{\pi \gamma^{-1} s} }  + { \pi^2\gamma^{-2}e^{-cs} \cos{\pi\gamma^{-1}s}
    \over \sin^2{\pi \gamma^{-1} s} }, \quad 0<s< \gamma.
\end{align*}
Without loss of generality, assume that $(\gamma_i, c_i)$,
$i=1,\ldots,k$, satisfy a lexicographical ordering, that is,
$\gamma_i\leq \gamma_{i+1}$ and $c_i<c_{i+1}$ when
$\gamma_i=\gamma_{i+1}$. Suppose that
$\gamma_1=\cdots=\gamma_l<\gamma_{l+1}\leq \cdots \leq \gamma_k$ and
$c_1<\cdots<c_l$. Then, applying the two-sided laplace
transformation to (\ref{Lemma equality}), we have that  for all
$0<s<\gamma_1$,
\begin{align}\label{Lemma transformed eqn}
    \sum_{i=1}^k d_{i0} F_0(s;\gamma_i,c_i) + \sum_{i=1}^k d_{i1} F_1(s;\gamma_i,c_i)=0.
\end{align}
Since the numerator of the left-hand side of (\ref{Lemma transformed
eqn}) is an analytic function on $\mathbb{R}$, (\ref{Lemma
transformed eqn}) is still valid for all $s\in\mathbb{R}\backslash
D$, where $D=\{ s: s=\gamma_i m, 1\leq i \leq k, m\in\mathbb{Z} \}$.
Multiplying $\sin^2{\pi\gamma_1^{-1}s}$ to both the sides of
(\ref{Lemma transformed eqn}), we attain
\begin{align*}
    & \sin{ \pi\gamma_1^{-1}s} \sum_{i=1}^l \left\{ d_{i0} \pi \gamma_1^{-1} e^{-c_i s}
    + d_{i1} \pi \gamma_1^{-1} c_i e^{-c_i s} \right\}  \\
    &+ \sin^2{ \pi\gamma_1^{-1}s} \sum_{i=l+1}^k \left\{ d_{i0} {\pi \gamma_i^{-1} e^{-c_i s} \over \sin{\pi \gamma_i^{-1} s} }
    + d_{i1} {\pi \gamma_i^{-1} c_i e^{-c_i s} \over \sin{ \pi\gamma_i^{-1}s} }
    \right\} \\
    &+ \cos{\pi\gamma_1^{-1}s} \sum_{i=1}^l d_{i1} \pi^2 \gamma_1^{-2} e^{-c_i s}
    + \sin^2{\pi\gamma_1^{-1}s} \sum_{i=l+1}^k d_{i1} {\pi^2 \gamma_i^{-2} e^{-c_i s} \cos{\pi\gamma_i^{-1}s} \over \sin^2{\pi\gamma_i^{-1}s} }
    =0.
\end{align*}
Then, if we set $\mathbb{N}_1 = \{ n\in\mathbb{N} : \gamma_1 n \ne
\gamma_i m\ {\rm for\ all}\ l< i\leq k, m\in\mathbb{N} \}$, for any
fixed $n\in\mathbb{N}_1$, letting $s\to \gamma_1 n$ through the
values in $\mathbb{R}\backslash D$, we can have
\begin{align}\label{Lemma equality2}
    \sum_{i=1}^l d_{i1}  e^{-c_i \gamma_1 n}=0.
\end{align}
Since (\ref{Lemma equality2}) holds for all $n\in\mathbb{N}_1$,
multiplying $e^{c_1 \gamma_1 n}$ to both the sides of (\ref{Lemma
equality2}) and letting $n\to\infty$ through the values in
$\mathbb{N}_1$, we get $d_{11}=0$. Similarly, it can be seen that
$d_{21}=\cdots=d_{l1}=0$. Meanwhile, multiplying
$\sin{\pi\gamma_1^{-1}s}$ to both the sides to (\ref{Lemma
transformed eqn})  and letting $s\to \gamma_1 n$, we can have $
\sum_{i=1}^l d_{i0} e^{-c_i \gamma_1 n}=0$ for any
$n\in\mathbb{N}_1$, and henceforth, $d_{10}=\cdots=d_{l0}=0$.
Continuing the above process, one can finally establish the
lemma.\hfill \fbox{}

\begin{remark}
    Lemma~\ref{Linear independence Lemma} actually entails the identifiability of
  logistic mixture distributions (cf. \cite{YakowitzSpragins1968} and \cite{Sussmann1992}).
  \cite{HwangDing1997} also proved the linear independence of
  logistic distributions  and their density functions to deal with the
  identifiability problem
  in artificial neural networks. However, their results do
  not directly imply Lemma~\ref{Linear independence Lemma}.
  Our proof is simpler and is based on Theorem~2 of \cite{Teicher1963}.
\end{remark}

\section*{Acknowledgements}
The authors thank the referees for their careful reading and
valuable comments. This work was supported by the National Research
Foundation of Korea(NRF) grant funded by the Korea government(MEST)
(No. 2012R1A2A2A01046092).

\end{document}